\documentclass[a4paper,12pt]{amsart}
\usepackage{amssymb,amscd}

 \calclayout

\theoremstyle{plain}
\newtheorem{theo}{Theorem}[section]
\newtheorem{lemm}[theo]{Lemma}
\newtheorem{coro}[theo]{Corollary}

\theoremstyle{remark}
\newtheorem*{rema}{Remark}

\newcommand{\field}[1]{\mathbb{#1}}
\newcommand{\C}{\field{C}}
\newcommand{\Z}{\field{Z}}

\newcommand{\side}[1]{{}^{#1}\hskip-.075em}

\DeclareMathOperator{\Hom}{Hom} \DeclareMathOperator{\Irr}{Irr}
\DeclareMathOperator{\Vect}{Vect} \DeclareMathOperator{\ind}{ind}
\DeclareMathOperator{\res}{res}

\begin{document}

\title{Equivariant $K$-groups of spheres with actions of involutions}

\author[Jin-Hwan Cho]{Jin-Hwan Cho$^\dagger$}
\address{School of Mathematics, Korea Institute for Advanced
Study \\ 207-43 Cheongryangri-dong, Dongdaemun-gu \\ Seoul
130-012, Republic of Korea}
\email{chofchof@kias.re.kr}

\author{Mikiya Masuda}
\address{Department of Mathematics, Osaka City University \\
3-3-138 Sugimoto, Sumiyoshi-ku \\ Osaka 558, Japan}
\email{masuda@sci.osaka-cu.ac.jp}

\thanks{The manuscript was written based on the authors' talk
at the 2nd Japan-Korea conference on Transformation Groups on
August 8, 1999 at Okayama University of Science}

\thanks{$^\dagger$The first author was supported by postdoctoral fellowships
program from Korea Science \& Engineering Foundation (KOSEF)}

\date{May 29, 2000}



\begin{abstract}
We calculate the $R(G)$-algebra structure on the reduced
equivariant $K$-groups of two-dimensional spheres on which a
compact Lie group $G$ acts as involutions. In particular, the
reduced equivariant $K$-groups are trivial if $G$ is abelian,
which shows that the previous Y. Yang's calculation
in~\cite{Yan95} is not true.
\end{abstract}

\maketitle

\section{Introduction}

Let $G$ be a compact Lie group. As one of (stable) equivariant
cohomology theories, $RO(G)$-graded cohomology theory has been
considered a full equivariant generalization of ordinary singular
cohomology, which is built in the interplay relating the Burnside
ring, the real character ring $RO(G)$, and $G$-homotopy theory
(see~\cite{LMM81} for more details). The main difference from the
other equivariant cohomology theories is that the coefficient
ring is indexed by all virtual real representations instead of
integers, and equivariant $K$-theory is one of the adequate known
examples.

In 1995, Y. Yang~\cite{Yan95} proved that, for any compact Lie
group $G$, the coefficient groups $\widetilde
K_G^*(S^0)\equiv\widetilde K_G(S^*)$ of $RO(G)$-graded equivariant
$K$-theory can only have $2$-torsion. He calculated the reduced
equivariant $K$-groups of $S^1$ and $S^2$ with a finite cyclic
group $G$ acting as involutions, and deduced his main result using
the Bott periodicity theorem and J.~E. McClure's
results~\cite{McC86}.

More precisely, if a compact $G$-space $X$ has a base point $*$
fixed by the $G$-action, then the reduced equivariant $K$-group
$\widetilde K_G(X)$ is defined to be the kernel of the restriction
homomorphism $K_G(X)\to K_G(*)$ induced from the inclusion map
$*\hookrightarrow X$. In fact, $K_G(X)$ and $\widetilde K_G(X)$
are algebras over the complex character ring $R(G)$ (although
there is no identity element in $\widetilde K_G(X)$). Let
$\lambda\colon G\to O(1)=\{\pm1\}$ be a surjective homomorphism
regarded as a one-dimensional real representation of $G$. Denote
by $1$ the trivial one-dimensional real representation of $G$.
Actually Y. Yang calculated $\widetilde K_G(S^\lambda)$ and
$\widetilde K_G(S^{1\oplus\lambda})$ for finite cyclic groups
$G$, where $S^\lambda$ and $S^{1\oplus\lambda}$ denote the
one-point compactifications of the real representations $\lambda$
and $1\oplus\lambda$, respectively.

The authors proved recently in~\cite[Theorem~10.1]{CKMS99} that
$\widetilde K_G(S^\lambda)$ is isomorphic to the ideal in $R(G)$
generated by $(1-\lambda)\otimes\C$, which extends Y. Yang's
result~\cite[Theorem A]{Yan95} for $G$ finite cyclic to any
compact Lie group. In this paper we apply the same technique to
calculate the $R(G)$-algebra structure of $\widetilde
K_G(S^{1\oplus\lambda})$ when $G$ is a compact Lie group. In
particular, we will prove that $\widetilde
K_G(S^{1\oplus\lambda})$ is trivial for any compact abelian Lie
group $G$, which shows that Y. Yang's result~\cite[Theorem
B]{Yan95} is not true.

In the following our main results are stated. Denote by $H$ the
kernel of the real representation $\lambda$. Given a character
$\chi$ of $H$ and $g\in G$, a new character $\side{g}\chi$ of $H$
is defined by $\side{g}\chi(h)=\chi(g^{-1}hg)$ for $h\in H$.
Choose and fix an element $b\in G\setminus H$. Since a character
is a class function and $G/H$ is of order two,
$\side{b}\chi=\side{g}\chi$ for all $g\in G\setminus H$ so that
$\side{b}\chi$ is independent of the choice of the element $b\in
G\setminus H$. Note that $R(H)$ has a canonical $R(G)$-module
structure given by $\varphi\cdot\chi=\res_H\varphi\otimes\chi$
for $\varphi\in R(G)$ and $\chi\in R(H)$.

\begin{theo} \label{theo:main_theorem}
Let $G$ be a compact Lie group and let $\lambda\colon G\to O(1)$
be a surjective homomorphism. Denote by $H$ the kernel of
$\lambda$ and choose an element $b\in G\setminus H$. Then
$\widetilde K_G(S^{1\oplus\lambda})$ is isomorphic to the
$R(G)$-module consisting of the elements $\chi-\side{b}\chi$ in
$R(H)$ for all characters $\chi$ of $H$. Moreover, the ring
structure on $\widetilde K_G(S^{1\oplus\lambda})$ is given by
$\alpha\beta=0$ for all elements $\alpha,\beta\in\widetilde
K_G(S^{1\oplus\lambda})$.
\end{theo}

\begin{coro} \label{coro:main_corollary}
In particular, $\widetilde K_G(S^{1\oplus\lambda})$ is trivial if
there exists an element in $G\setminus H$ commuting with all
elements in $H$ (for example, $G$ is abelian).
\end{coro}

\section{Complex $\Z_2$-vector bundles over $S^{1\oplus\lambda}$}

We will begin by considering the structure of complex
$\Z_2$-vector bundles over $S^{1\oplus\lambda}$. In this case
$\lambda\colon\Z_2\to O(1)=\{\pm1\}$ becomes an isomorphism.

\begin{lemm} \label{lemm:decomposition}
Every complex $\Z_2$-vector bundle over $S^{1\oplus\lambda}$
decomposes into the Whitney sum of $\Z_2$-invariant sub-line
bundles.
\end{lemm}

\begin{proof}
The set of $\Z_2$-fixed points in $S^{1\oplus\lambda}$ constitutes
a circle, denoted by $S^1$, containing $0$ and $*$. Then $S^1$
divides the sphere $S^{1\oplus\lambda}$ into two hemispheres, say
$H_1$ and $H_2$. Let $E$ be a complex $\Z_2$-vector bundle over
$S^{1\oplus\lambda}$. Note that the restriction of $E$ to the
subspace $S^1\subset S^{1\oplus\lambda}$ decomposes into the
Whitney sum of $\Z_2$-invariant sub-line bundles. Choose a
$\Z_2$-invariant sub-line bundle, say $F$, over $S^1$. Then it is
always possible to extend $F$ to a non-equivariant sub-line
bundle over $H_1$, since the restriction of $E$ to $H_1$ is
non-equivariantly trivial and the set of one-dimensional
subspaces of the fiber is homeomorphic to the Grassmann manifold
$G_\C(k,1)\cong\C P^{k-1}$, the fundamental group of which is
trivial. We now extend it over the other hemisphere $H_2$ using
the $\Z_2$-action on $E$ to get a resulting $\Z_2$-invariant
sub-line bundle of $E$.
\end{proof}

\begin{lemm} \label{lemm:Z_2-triviality}
Every complex $\Z_2$-vector bundle over $S^{1\oplus\lambda}$ is
trivial.
\end{lemm}

\begin{proof}
It suffices to show that every complex $\Z_2$-line bundle $E$ over
$S^{1\oplus\lambda}$ is trivial by Lemma~\ref{lemm:decomposition}.
Choose two $\Z_2$-invariant hemispheres, denoted by $H_0$ and
$H_*$ of $S^{1\oplus\lambda}$ containing $0$ and $*$,
respectively, such that $H_0\cup H_*=S^{1\oplus\lambda}$ and
$H_0\cap H_*$ is a $\Z_2$-invariant circle. Since $H_0$
(resp.~$H_*$) is equivariantly contractible to $0$ (resp.~$*$),
$E$ restricted to $H_0$ (resp.~$H_*$) is equivariantly isomorphic
to the product bundle $H_0\times E_0$ (resp.~$H_*\times E_*$)
where $E_0$ (resp.~$E_*$) denotes the fiber of $E$ at $0$
(resp.~$*$). Then $E$ induces a $\Z_2$-equivariant clutching map
of $H_0\cap H_*$ into the set $\Hom(E_0,E_*)$ of linear maps
between the two fibers $E_0$ and $E_*$. Note that $E_0$ and $E_*$
are isomorphic as complex representations of $\Z_2$ since $0$ and
$*$ are connected by $\Z_2$-fixed points, and thus the induced
$\Z_2$-action on $\Hom(E_0,E_*)$ is trivial. Since $\Z_2$ acts on
the circle $H_0\cap H_*$ as a reflection, the clutching map is
equivariantly null-homotopic, that is, $E$ is equivariantly
trivial.
\end{proof}

\section{Induced vector bundles}

In this section we review the notion of induced vector bundles
defined in a way similar to the definition of induced
representations.

Let $G$ be a compact Lie group and let $X$ be a $G$-space. The
set of isomorphism classes of (real or complex) $G$-vector
bundles over $X$, usually denoted by $\Vect_G(X)$, is a semi-group
under the Whitney sum operation. For a closed subgroup $H$ of
$G$, the inclusion map $H\hookrightarrow G$ induces a semi-group
homomorphism $\res_H\colon\Vect_G(X)\to\Vect_H(X)$ called the
\emph{restriction} homomorphism. On the other hand, there is
another homomorphism $\ind_H^G\colon\Vect_H(X)\to\Vect_G(X)$
called the \emph{induction} homomorphism (see for
instance~\cite{ChSu85} or~\cite{MaPe85}). In addition, if $H$ has
finite index in $G$, then the induction homomorphism can be
defined explicitly, which we will explain below.

Let $E$ be an $H$-vector bundle over $X$. Choose a set of
representatives $\{t_0,t_1,\dotsc,t_l\}$ of $G/H$. The identity
element of $G$ will always be selected as $t_0$. Define a group
isomorphism $\imath_g\colon gHg^{-1}\to H$ by
$\imath_g(ghg^{-1})=h$ for $h\in H$ and $g\in G$.

We first construct a $t_iHt_i^{-1}$-vector bundle $t_iE$ over $X$
for each $0\leq i\leq l$. Consider the following pull-back diagram
\[
\begin{CD}
(t_i^{-1})^*E @>>> E \\ @VVV @VVV \\ X @>{t_i^{-1}}>> X,
\end{CD}
\]
where $t_i^{-1}$ denotes the action of $t_i^{-1}$ on $X$. Since
the action of $t_i^{-1}$ on $X$ is $\imath_{t_i}$-equivariant,
i.e., $t_i^{-1}(hx)=\imath_{t_i}(h)t_i^{-1}(x)$ for $h\in
t_iHt_i^{-1}$ and $x\in X$, the pull-back bundle $(t_i^{-1})^*E$
has a canonical $t_iHt_i^{-1}$-action and it becomes a
$t_iHt_i^{-1}$-vector bundle over $X$, simply denoted by $t_iE$.
Note that every element in the fiber of $t_iE$ at $x\in X$ is
represented by $t_iv$ for some $v$ in the fiber of $E$ at
$t_i^{-1}x$. We now define $\ind_H^G E$ by the Whitney sum
(without action) of all the bundles $t_iE$ and then give a
$G$-action on it as follows.

Given $g\in G$, each element $gt_i$ has a unique form $t_{i'}h_i$
for some representative $t_{i'}\in\{t_0,\dotsc,t_l\}$ and $h_i\in
H$. So we define an action of $g$ on $\ind_H^G E$ by
\[
g\cdot\Bigl(\sum t_i v_i\Bigr) = \sum t_{i'} (h_i v_i).
\]
If $\sum t_i v_i$ is an element in the fiber of $\ind_K^G E$ at
$x\in X$, then each $v_i$ is contained in the fiber of $E$ at
$t_i^{-1}x$ and $h_iv_i$ is in the fiber of $E$ at $t_{i'}^{-1}gx$
since $h_it_i^{-1}=t_{i'}^{-1}g$. Thus the action of $g$ gives a
linear map between the two fibers at $x$ and $gx$ for all $x\in
X$. It is easy to see that the definition gives actually an
action of $G$ on $\ind_K^G E$. Moreover the action is compatible
with the $t_iHt_i^{-1}$-action already defined on $t_iE$.

We next show that the induced bundle is independent of the choice
of the representatives $\{t_0,\dotsc,t_l\}$ of $G/H$. Choose
another set of representatives $\{s_0=e,s_1,\dotsc,s_l\}$ of
$G/H$. We may assume that each $s_i^{-1}t_i$ is contained in $H$
by arranging the order of representatives if necessary, and that
every element in the induced bundle constructed with
$\{s_0,\dotsc,s_l\}$ is written by $\sum s_iv_i$. It is easy to
check that the bundle map defined by $\sum t_iv_i \mapsto \sum
s_i(s_i^{-1}t_i)v_i$ gives an equivariant isomorphism between the
two induced bundles.

\begin{lemm} \label{lemm:uniqueness_of_induced_bundles}
If two $H$-vector bundles $E$ and $E'$ over $X$ are isomorphic,
then so are the induced bundles $\ind_H^G E$ and $\ind_H^G E'$.
The converse is not true in general.
\end{lemm}

\begin{proof}
Given an $H$-vector bundle isomorphism $\Phi_0\colon E\to E'$,
the map
\[
\Phi\colon\ind_H^G E\to\ind_H^G E'
\]
sending $\sum t_iv_i$ to $\sum t_i\Phi_0(v_i)$ gives a $G$-vector
bundle isomorphism.
\end{proof}

\begin{lemm} \label{lemm:triviality_of_induced_bundles}
Let $W$ be a representation of $H$. Then the induced bundle of the
product bundle $X\times W\to X$ is isomorphic to the product
bundle $X\times \ind_H^G W\to X$.
\end{lemm}

\begin{proof}
The construction of induced bundles is the same as that of induced
representations.
\end{proof}

\begin{lemm} \label{lemm:tesoring_property}
Let $L$ be a $G$-vector bundle and let $E$ be an $H$-vector
bundle over the same base space $X$. Then $L\otimes\ind_H^GE$ is
isomorphic to $\ind_H^G(\res_HL\otimes E)$.
\end{lemm}

\begin{proof}
For $l\in L$ and $\sum t_iv_i\in\ind_H^GE$, the map
\[
\Phi\colon L\otimes\ind_H^GE\to\ind_K^G(\res_HL\otimes E)
\]
sending $l\otimes\sum t_iv_i$ to $\sum t_i(t_i^{-1}l\otimes v_i)$
gives a (non-equivariant) vector bundle map. Given $g\in G$ we
have $gt_i=t_{i'}h_i$ for some representative $t_{i'}$ and
$h_i\in H$. Then the equalities
\begin{align*}
\Phi\bigl(g\bigl(l\otimes\sum t_iv_i\bigr)\bigr) &=
\Phi\bigl(gl\otimes\sum t_{i'}(h_iv_i)\bigr)
= \sum t_{i'}\bigl(t_{i'}^{-1}gl\otimes h_iv_i\bigr) \\
&= \sum t_{i'}h_i\bigl(t_i^{-1}l\otimes v_i\bigr)
= \sum gt_i\bigl(t_i^{-1}l\otimes v_i\bigr) \\
&= g\Phi\bigl(l\otimes\sum t_iv_i\bigr)
\end{align*}
imply that $\Phi$ is $G$-equivariant. On the other hand, the
inverse bundle map is given by the map sending $\sum
t_i(l_i\otimes v_i)$ to $\sum(t_il_i\otimes t_iv_i)$ for $l_i\in
L$ and $v_i\in E$.
\end{proof}

\section{Decomposition of equivariant vector bundles}

In this section we rephrase the relevant material
from~\cite[Section~2]{CKMS99} to decompose equivariant vector
bundles for the readers' convenience.

Let $G$ be a compact Lie group and $H$ a closed normal subgroup
of $G$. Given a character $\chi$ of $H$ and $g\in G$, a new
character $\side{g}\chi$ of $H$ is defined by
$\side{g}\chi(h)=\chi(g^{-1}hg)$ for $h\in H$. This defines an
action of $G$ on the set $\Irr(H)$ of irreducible characters of
$H$. Since a character is a class function, $H$ acts trivially on
$\Irr(H)$. Therefore, the isotropy subgroup of $G$ at
$\chi\in\Irr(H)$, denoted by $G_\chi$, contains $H$. We choose a
representative from each $G$-orbit in $\Irr(H)$ and denote by
$\Irr(H)/G$ the set of those representatives.

Let $X$ be a connected $G$-space on which $H$ acts trivially.
Then all the fibers of a $G$-vector bundle $E$ over $X$ are
isomorphic as representations of $H$. We call the unique (up to
isomorphism) representation of $H$ the \emph{fiber
$H$-representation} of $E$.

As is well-known, $E$ decomposes according to irreducible
representations of $H$. For $\chi\in\Irr(H)$, we denote by
$E(\chi)$ the $\chi$-isotypical component of $E$, that is, the
largest $H$-subbundle of $E$ with a multiple of $\chi$ as the
character of the fiber $H$-representation. Note that $gE(\chi)$,
that is $E(\chi)$ mapped by $g\in G$, is
$\side{g}\chi$-isotypical component of $E$. This means that
$E(\chi)$ is actually a $G_\chi$-vector bundle and that
$\bigoplus_{\chi'\in G(\chi)}E(\chi')$, where $G(\chi)$ denotes
the $G$-orbit or $\chi$, is a $G$-subbundle of $E$. Since
$\bigoplus_{\chi'\in G(\chi)}E(\chi')$ is nothing but the induced
bundle $\ind_{G_\chi}^GE(\chi)$, we have the following
decomposition
\begin{equation}
E = \bigoplus_{\chi\in\Irr(H)/G} \ind_{G_\chi}^G E(\chi) \tag{*}
\label{eq:bundle_decomposition}
\end{equation}
as $G$-vector bundles.

\begin{lemm} \label{lemm:orbit_isomorphism}
Two $G$-vector bundles $E$ and $E'$ are isomorphic if and only if
$E(\chi)$ and $E'(\chi)$ are isomorphic as $G_\chi$-vector
bundles for each $\chi\in\Irr(H)/G$. In particular, $E$ is
trivial if and only if $E(\chi)$ is trivial for each
$\chi\in\Irr(H)/G$.
\end{lemm}

\begin{proof}
The necessity is obvious since a $G$-vector bundle isomorphism
$E\to E'$ restricts to a $G_\chi$-vector bundle isomorphism
$E(\chi)\to E'(\chi)$, and the sufficiency follows from the fact
that $\ind_{G_\chi}^G$ is functorial.
\end{proof}

\begin{coro} \label{coro:orbit_isomorphism}
Two $G$-vector bundles $E$ and $E'$ are stably isomorphic if and
only if $E(\chi)$ and $E'(\chi)$ are stably isomorphic for each
$\chi\in\Irr(H)/G$. \qed
\end{coro}

The observation above can be restated in $K$-theory as follows.
Denote by $K_{G_\chi}(X,\chi)$ the subgroup of $K_{G_\chi}(X)$
generated by $G_\chi$-vector bundles over $X$ with a multiple of
$\chi$ as the character of fiber $H$-representations. The reduced
version $\widetilde K_{G_\chi}(X,\chi)$ can be defined
accordingly. Then the map sending $E$ to
$\prod_{\chi\in\Irr(H)/G}E(\chi)$ gives group isomorphisms
\[
K_G(X) \to \prod_{\chi\in\Irr(H)/G} K_{G_\chi}(X,\chi)
\quad\text{and}\quad \widetilde K_G(X) \to
\prod_{\chi\in\Irr(H)/G} \widetilde K_{G_\chi}(X,\chi)
\]
by Corollary~\ref{coro:orbit_isomorphism}.

\begin{lemm} \label{lemm:HOM_argument}
If there is a complex $G_\chi$-vector bundle $L$ over $X$ with
$\chi$ as the character of the fiber $H$-representation, then the
map sending $E$ to $\Hom_H(L,E)$ gives group isomorphisms
\[
K_{G_\chi}(X,\chi)\to K_{G_\chi/H}(X) \quad\text{and}\quad
\widetilde K_{G_\chi}(X,\chi)\to \widetilde K_{G_\chi/H}(X).
\]
\end{lemm}

\begin{proof}
It is easy to check that the map
\[
K_{G_\chi/H}(X)\to K_{G_\chi}(X,\chi)
\]
sending $F$ to $L\otimes F$ gives the inverse of the map in the
lemma.
\end{proof}

\begin{rema}
The lemma above does not hold in the real category in general,
but it does if $\chi$ is the character of a real irreducible
representation of $H$ with the endomorphism algebra isomorphic to
the set of real numbers.
\end{rema}

\section{Proof of Theorem~\ref{theo:main_theorem}}

We now return to our original setting to prove the main result.
Hereafter we omit the adjective ``complex'' for complex vector
bundles and complex representations since we work in the complex
category. At first consider the additive structure on $\widetilde
K_G(S^{1\oplus\lambda})$.

Let $G$ be a compact Lie group. Denote by $H$ the kernel of the
surjective homomorphism $\lambda\colon G\to O(1)=\{\pm1\}$.
Choose and fix an element $b\in G\setminus H$. Since $G_\chi$
contains $H$ and $G/H$ is of order two, $G_\chi$ is either $H$ or
$G$ for each irreducible character $\chi\in\Irr(H)/G$. Note that,
in each case, the condition of Lemma~\ref{lemm:HOM_argument}
holds, since there exists a $G_\chi$-extension of
$\chi$~\cite[Proposition~4.2]{CKS99} and it gives a trivial
$G_\chi$-vector bundle with $\chi$ as the character of the fiber
$H$-representation. Therefore, according to the arguments in the
previous section, we have a group isomorphism
\begin{equation}
\widetilde K_G(S^{1\oplus\lambda}) \cong
\prod_{\substack{\chi\in\Irr(H)/G\\\text{s.t.~$G_\chi=H$}}}
\widetilde K(S^2) \times
\prod_{\substack{\chi\in\Irr(H)/G\\\text{s.t.~$G_\chi=G$}}}
\widetilde K_{G/H}(S^{1\oplus\lambda}). \tag{**}
\label{eq:the_structure_of K_group}
\end{equation}
It is well-known that $\widetilde K(S^2)$ is infinite cyclic with
generator $\zeta-1$, where $\zeta$ and $1$, respectively, denote
the dual bundle of the canonical line bundle and the trivial line
bundle over $S^2\cong\C P^1$ (see for
instance~\cite[Theorem~2.3.14]{Ati67}). Moreover, $\widetilde
K_{G/H}(S^{1\oplus\lambda})$ is trivial since $G/H\cong\Z_2$ and
every $\Z_2$-vector bundle over $S^{1\oplus\lambda}$ is trivial
by~Lemma~\ref{lemm:Z_2-triviality}. Therefore the
decomposition~(\ref{eq:bundle_decomposition})
and~Lemma~\ref{lemm:HOM_argument} in the previous section imply
the following lemma.

\begin{lemm} \label{lemm:group_structure}
$\widetilde K_G(S^{1\oplus\lambda})$ is a free abelian group
generated by $\ind_H^G\bigl(\chi\otimes(\zeta-1)\bigr)$ for each
$\chi\in\Irr(H)/G$ such that $G_\chi=H$. Here, the same notation
$\chi$ is used for the product $H$-vector bundle over
$S^{1\oplus\lambda}$ with $\chi$ as the character of the fiber
representation. \qed
\end{lemm}

\begin{proof}[Proof of Theorem~\ref{theo:main_theorem}]
We now consider the map
\[
\Psi\colon\widetilde K_G(S^{1\oplus\lambda})\to R(H)
\]
sending each generator $\ind_H^G\bigl(\chi\otimes(\zeta-1)\bigr)$
of $\widetilde K_G(S^{1\oplus\lambda})$ to $\chi-\side{b}\chi\in
R(H)$. Since both $\widetilde K_G(S^{1\oplus\lambda})$ and $R(H)$
are free abelian groups, $\Psi$ is a well-defined group
homomorphism. Moreover it is injective, since the set of the
elements $\chi-\side{b}\chi$ for all $\chi\in\Irr(H)/G$ such that
$G_\chi=H$ can be extended to an additive basis of $R(H)$.

For each $\chi\in\Irr(H)$, either $G_\chi=H$ or $G$. In case that
$G_\chi=G$, we have $\chi-\side{b}\chi=0$. Thus the image of
$\Psi$ is generated by the elements $\chi-\side{b}\chi$ for all
$\chi\in\Irr(H)$, that is, the $R(G)$-submodule of $R(H)$
consisting of the elements $\chi-\side{b}\chi$ for all characters
$\chi$ of $H$.

Given a character $\varphi$ of $G$, as in
Lemma~\ref{lemm:group_structure}, we use the same notation
$\varphi$  for the product $G$-vector bundle over
$S^{1\oplus\lambda}$ with $\varphi$ as the character of the fiber
representation. Then Lemma~\ref{lemm:tesoring_property} implies
that
\[
\varphi\otimes\ind_H^G\bigl(\chi\otimes(\zeta-1)\bigr) =
\ind_H^G\bigl(\res_H\varphi\otimes\chi\otimes(\zeta-1)\bigr).
\]
Since $\side{b}(\res_H\varphi)=\res_H\varphi$, we have the
equalities
\begin{align*}
\Psi\bigl(\varphi\otimes\ind_H^G(\chi\otimes(\zeta-1))\bigr)&=\res_H\varphi\otimes\chi-\side{b}(\res_H\varphi)\otimes\side{b}\chi\\
&=\res_H\varphi\otimes(\chi-\side{b}\chi)\\
&=\varphi\cdot\Psi\bigl(\ind_H^G(\chi\otimes(\zeta-1))\bigr)
\end{align*}
showing that $\Psi$ is an $R(G)$-module homomorphism.

It remains to show the ring structure on $\widetilde
K_G(S^{1\oplus\lambda})$. It suffices to show that the tensor
product of any two generators in $\widetilde
K_G(S^{1\oplus\lambda})$ is zero. Note that, given an induced
bundle $\ind_H^G(\chi\otimes\zeta)$, the image of
$\chi\otimes\zeta$ mapped by the $b$-action, that is the
$\side{b}\chi$-isotypical component of
$\ind_H^G(\chi\otimes\zeta)$, is isomorphic to
$\side{b}\chi\otimes\zeta^*$ where $\zeta^*$ denotes the dual
bundle of $\zeta$. Indeed, the pull-back bundle
$(b^{-1})^*(\zeta)$ is isomorphic to $\zeta^*$, since the action
$b^{-1}\colon S^{1\oplus\lambda}\to S^{1\oplus\lambda}$ is a
reflection so that it induces the multiplication by $-1$ on the
second cohomology level. It follows that
\begin{align*}
\res_H\ind_H^G\bigl(\chi\otimes(\zeta-1)\bigr)&=\chi\otimes(\zeta-1)+\side{b}\chi\otimes(\zeta^*-1)\\
&=\chi\otimes(\zeta-1)-\side{b}\chi\otimes(\zeta-1)\\
&=(\chi-\side{b}\chi)\otimes(\zeta-1),
\end{align*}
since $\zeta^*-1=-(\zeta-1)$ in $\widetilde K(S^2)$. Therefore,
by~Lemma~\ref{lemm:tesoring_property}, we have the equalities
\begin{align*}
&\ind_H^G\bigl(\chi\otimes(\zeta-1)\bigr)\otimes\ind_H^G\bigl(\eta\otimes(\zeta-1)\bigr)\\
&\qquad\qquad=\ind_H^G\bigl(\res_H\ind_H^G(\chi\otimes(\zeta-1))\otimes\eta\otimes(\zeta-1)\bigr)\\
&\qquad\qquad=\ind_H^G\bigl(\bigl((\chi-\side{b}\chi)\otimes\eta\otimes(\zeta-1)^2\bigr)\\
&\qquad\qquad= 0,
\end{align*}
since $(\zeta-1)^2=0$ in $\widetilde K(S^2)$
(see~\cite[Theorem~2.3.14]{Ati67}).
\end{proof}

\begin{proof}[Proof of Corollary~\ref{coro:main_corollary}]
By assumption there exists an element $b\in G\setminus H$ such
that $bh=hb$ for all $h\in H$. It follows that
$\side{b}\chi(h)=\chi(b^{-1}hb)=\chi(h)$ for any character $\chi$
of $H$, which completes the proof.
\end{proof}

\providecommand{\bysame}{\leavevmode\hbox to3em{\hrulefill}\thinspace}

\end{document}